\documentclass[conference]{IEEEtran}
\IEEEoverridecommandlockouts

\ifCLASSINFOpdf
\else
\fi

\usepackage{graphicx}
\usepackage[cmex10]{amsmath}
\usepackage{amssymb}
\usepackage{epstopdf}
\graphicspath{{./fig/}}
\usepackage{pgfplots}
\usepackage{enumitem}
\usepackage{array}
\IEEEoverridecommandlockouts

\newcommand{\ve}{\mathbf}
\newcommand{\m}{\mathbf}

\newcommand{\mf}[1]{\mathbf{\tilde{\mathbf{#1}}}} 
\newcommand{\vef}[1]{\mathbf{\tilde{\mathbf{#1}}}} 

\begin{document}
%
\title{On the Log-Likelihood Ratio Evaluation of CWCU Linear and Widely Linear MMSE Data Estimators}

\author{\IEEEauthorblockN{Oliver~Lang and Mario~Huemer \thanks{Copyright 2001 SS\&C. Published in the Proceedings of the \emph{Asilomar Conderence on Signals, Systems, and Computers}, November 6-9th, 2016, Pacific Grove, CA, USA. This work has been supported by the Austrian Science Fund (FWF): I683-N13.}}
\IEEEauthorblockA{Johannes Kepler University\\
	Institute of Signal Processing\\
	4040 Linz}
\and
\IEEEauthorblockN{Christian~Hofbauer}
\IEEEauthorblockA{Linz Center of Mechatronics GmbH\\
	4040 Linz}}


%


\maketitle

\begin{abstract}
In soft decoding, log-likelihood ratios (LLRs) are calculated from estimated data symbols. Data symbols from proper constellation diagrams such as QPSK are often estimated using the linear minimum mean square error (LMMSE) estimator. We prove that the recently introduced component-wise conditionally unbiased (CWCU) LMMSE estimator results in the very same LLRs as the LMMSE estimator for typical model assumptions. For improper constellation diagrams such as 8-QAM, we show that the widely linear versions of the LMMSE and the CWCU LMMSE estimator also yield identical LLRs. In that case, the CWCU estimator allows to reduce the complexity of the LLR determination.

\end{abstract}
%

%
\IEEEpeerreviewmaketitle

\section{Introduction} \label{Sec:Introduction}
%
%
%
%
The task of estimating a parameter vector $\ve{x}\in\mathbb{C}^{n \times 1}$ out of a measurement vector $\ve{y}\in\mathbb{C}^{m \times 1}$ with $m \geq n$ can be treated in the classical sense or in the Bayesian sense. Classical and Bayesian estimation not only differ in terms of the incorporation of prior knowledge, but also in terms of the unbiased properties. While a \underline{c}lassical estimator $\hat{\ve{x}}_{\text{C}}$ has to fulfill
\begin{equation}
	E_\ve{y}[\hat{\ve{x}}_{\text{C}}] = \ve{x} \hspace{0.4cm} \text{for all possible } \ve{x} \label{equ:cuLMMSE045}
\end{equation}
to be considered as unbiased, the frequently applied Bayesian linear minimum mean square error (\underline{L}MMSE) estimator only fulfills
\begin{equation}
	E_{\ve{y},\ve{x}}[\hat{\ve{x}}_{\text{L}} - \ve{x}] = E_{\ve{x}}\left[E_{\ve{y}|\ve{x}}\left[\hat{\ve{x}}_{\text{L}} - \ve{x}|\ve{x}\right] \right] = \ve{0}. \label{equ:cuLMMSE045a}
\end{equation}
This means $\hat{\ve{x}}_{\text{L}}$ is only "unbiased" when averaged over the probability density function (PDF) of $\ve{x}$, which is a much weaker constraint than \eqref{equ:cuLMMSE045}. However, the Bayesian approach allows the incorporation of prior knowledge. In \cite{Slock-Oct2005}--\cite{Lang_2015_Eurocast}, an interesting compromise between the stringent classical unbiased constraint and the weak Bayesian unbiased constraint has been investigated. 
There, component-wise conditionally unbiased (CWCU) Bayesian parameter estimators have been studied, which aim for achieving conditional unbiasedness for one parameter component at a time. Let $x_i$ be the $i^{th}$ element of $\ve{x}$ and $\hat{x}_i$ be an estimator of $x_i$, 
then the CWCU constraints are
\begin{equation}
	E_{\ve{y}|x_i}[\hat{x}_i|x_i] =  x_i, \label{equ:cuLMMSE047}  
\end{equation}
for all possible $x_i$ (and all $i=1,2,...,n$). The CWCU constraints are less stringent than the classical unbiased constraints in \eqref{equ:cuLMMSE045}, and it turns out that in many cases a CWCU estimator allows the incorporation of prior knowledge on the statistical properties of the parameter vector \cite{Huemer_2014_Asilomar}, \cite{Lang_2015_Eurocast}. In the following, we denote the linear estimator fulfilling the CWCU constraints and minimizing the Bayesian mean square error (BMSE) for $i = 1,2,\hdots,n$ as the CWCU LMMSE estimator. 
The CWCU LMMSE estimator is designed for proper measurement vectors. 
For the definition of propriety we refer to \cite{Schreier-2011}. 
A proper measurement vector could, e.g., arise when a data vector with proper symbols, such as for quadrature phase-shift keying (QPSK), is transmitted over a dispersive linear channel and disturbed by additive white Gaussian noise (AWGN). For this case the well-known LMMSE estimator is often used to estimate the transmitted symbols, followed by an evaluation of the log-likelihood ratios (LLRs). In Sec.~\ref{sec:Linear estimation of proper data} of this work it will be proven that the LLRs of the CWCU LMMSE estimates and the LMMSE estimates are identical even though the CWCU LMMSE estimator performs worse in terms of the BMSE. The second part of this paper focuses on improper symbol constellations such as 8 quadrature amplitude modulation (8-QAM). In such a scenario the widely LMMSE (WLMMSE) estimator is more appropriate for estimating the transmitted symbols. In Sec.~\ref{sec:Widely Linear Estimation of Improper Data} we prove that the CWCU WLMMSE estimator derived in \cite{CWCU_WLMMSE_Journal_Huemer_2016} again results in the same LLRs as the WLMMSE estimator while featuring a complexity advantage in deriving the LLR values. Finally, a simulation example is given in Sec.~\ref{sec:Propriety of the estimated data symbols} which illustrates the estimators' properties.

\section{LLR evaluation of proper symbols} \label{sec:Linear estimation of proper data}
In this section, the LLRs of proper symbols evaluated from the LMMSE estimates are compared with those determined from the CWCU LMMSE estimates. Let $\ve{x}$ and $\ve{y}$ be connected via the linear model $\ve{y} =\m{H}\ve{x} + \ve{n}$, 
where $\m{H}\in\mathbb{C}^{m\times n}$ is a known observation matrix, $\ve{x}$ has mean $E_{\ve{x}}[\ve{x}]$ and covariance matrix $\m{C}_{\ve{x}\ve{x}}=E_{\ve{x}}\left[(\ve{x}-E_{\ve{x}}[\ve{x}])(\ve{x}-E_{\ve{x}}[\ve{x}])^H\right]$ with $(\cdot )^H$ denoting the conjugate transposition, and $\ve{n}\in \mathbb{C}^{m \times 1}$ is a zero mean proper noise vector with covariance matrix $\m{C}_{\ve{n}\ve{n}}$ and independent of $\ve{x}$. Furthermore, let $\ve{h}_i\in\mathbb{C}^{m\times 1}$ be the $i^{th}$ column of $\m{H}$, $\bar{\m{H}}_i\in\mathbb{C}^{m\times (n-1)}$ the matrix resulting from $\m{H}$ by deleting $\ve{h}_i$, $x_i$ be the $i^{th}$ element of $\ve{x}$, and $\bar{\ve{x}}_i\in\mathbb{C}^{(n-1) \times 1}$ the vector resulting from $\ve{x}$ after deleting $x_i$. Then we can rewrite the linear model as
\begin{equation}
	\ve{y} = \ve{h}_i x_i + \bar{\m{H}}_i \bar{\ve{x}}_i + \ve{n}. \label{equ:CWCULMMSE015a} 
\end{equation}
Consider the general linear estimator $\hat{\ve{x}} = \m{E} \ve{y} , \hspace{0.2cm} \m{E}\in\mathbb{C}^{n\times m}$.
The $i^{th}$ component of this estimator is given by $\hat{x}_i = \ve{e}_i^H \ve{y}$, 
where $\ve{e}_i^H\in \mathbb{C}^{1 \times m}$ denotes the $i^{th}$ row of the estimator matrix $\m{E}$. Incorporating \eqref{equ:CWCULMMSE015a} yields
\begin{equation}
	\hat{x}_i = \underbrace{\ve{e}_i^H\ve{h}_i x_i}_{\text{Scaling}} + 
							 \underbrace{\ve{e}_i^H\bar{\m{H}}_i \bar{\ve{x}}_i}_{\text{IPI}} + 
						   \underbrace{\ve{e}_i^H\ve{n}}_{\text{Noise}}. \label{equ:CWCULMMSE015}
\end{equation}
In \eqref{equ:CWCULMMSE015}, we clearly see three effects, namely a scaling of the true parameter value, an inter-parameter interference (IPI) term, and a noise term. In communications, the noise term is usually Gaussian and the IPI term can usually approximately assumed to be Gaussian if $n$ is large enough due to central limit theorem arguments. From \eqref{equ:CWCULMMSE015}, the conditional mean of $\hat{x}_i$ becomes
\begin{equation}
	E_{\ve{y}|x_i}[\hat{x}_i|x_i] = \ve{e}_i^H\ve{h}_i x_i + 
					\ve{e}_i^H\bar{\m{H}}_i E_{\bar{\ve{x}}_i|x_i}[\bar{\ve{x}}_i|x_i]. \label{equ:CWCULMMSE016}
\end{equation} 
In the following, we assume statistically independent elements of $\ve{x}$ with zero mean, as usual in communications. Then, \eqref{equ:CWCULMMSE016} simplifies to
\begin{equation}
	E_{\ve{y}|x_i}[\hat{x}_i|x_i] = \ve{e}_i^H\ve{h}_i x_i = \alpha_i x_i. \label{equ:CWCULMMSE016a}
\end{equation}
The conditional variance of the general linear estimator is given by
\begin{align}
& \text{var}(\hat{x}_i|x_i)  \nonumber \\
 & \hspace{2mm}= E_{\ve{y}|x_i}\left[\left( \hat{x}_i - E_{\ve{y}|x_i}[\hat{x}_i|x_i]  \right) \left( \hat{x}_i - E_{\ve{y}|x_i}[\hat{x}_i|x_i]  \right)^H\middle|x_i \right].  \nonumber 
\end{align}
Inserting \eqref{equ:CWCULMMSE015}  and \eqref{equ:CWCULMMSE016a} into the previous equation yields
\begin{align}
\text{var}(\hat{x}_i|x_i) =& E_{\ve{y}|x_i}\left[( \ve{e}_i^H ( \bar{\m{H}}_i \bar{\ve{x}}_i + \ve{n}) ) ( \ve{e}_i^H ( \bar{\m{H}}_i \bar{\ve{x}}_i + \ve{n})  )^H|x_i \right] \nonumber \\
=&  \ve{e}_i^H ( \bar{\m{H}}_i \m{C}_{\bar{\ve{x}}_i \bar{\ve{x}}_i} \bar{\m{H}}_i^H + \m{C}_{\ve{n}\ve{n}}) \ve{e}_i. \label{equ:CWCULMMSE102} 
\end{align}
Note that the conditional variance in \eqref{equ:CWCULMMSE102} is independent of $x_i$. For a general estimator, the LLRs of any symbol constellation with equiprobable symbols can be written as \cite{Allpress_2004}
\begin{equation}
\Lambda(b_{ki}|\hat{x}_i) = \text{log} \frac{   \text{Pr}(b_{ki}=1|\hat{x}_i)}{   \text{Pr}(b_{ki}=0|\hat{x}_i)}= \text{log} \frac{\sum\limits_{q\in \mathit{S}(b_{ki}=1)}    p(\hat{x}_i|s^{(q)})}{\sum\limits_{q\in \mathit{S}(b_{ki}=0)}    p(\hat{x}_i|s^{(q)})}, \label{equ:LLR001}
\end{equation}
where $\hat{x}_i$ is the $i^{th}$ estimated symbol, $b_{ki}$ is the $k^{th}$ bit of the $i^{th}$ estimated symbol, $\mathit{S}(b_{ki}=1)$ and $\mathit{S}(b_{ki}=0)$ are the sets of symbol indices corresponding to $b_{ki}=1$ and $b_{ki}=0$, respectively, and $s^{(q)}$ is the $q^{th}$ symbol of such a set. 
In \eqref{equ:LLR001}, $ p(\hat{x}_i|s^{(q)})$ denotes the conditional PDF of the estimate $\hat{x}_i$ given that the actual symbol was $s^{(q)}$. Its Gaussian approximation is determined by the conditional mean and the conditional variance according to
\begin{equation}
 p(\hat{x}_i|s^{(q)}) = \frac{1}{\pi \text{var}(\hat{x}_i|s^{(q)})} e^{-\frac{1}{ \text{var}(\hat{x}_i|s^{(q)})} \left| \hat{x}_i -  E[\hat{x}_i|s^{(q)}] \right|^2}.  \label{equ:CWCULMMSE105}
\end{equation}
Together with \eqref{equ:LLR001}, the LLRs of any linear estimator can be evaluated by inserting the conditional mean and the conditional variance of the specific estimator. Such a specific estimator e.g., could be the LMMSE or the CWCU LMMSE estimators. We begin with the \underline{L}MMSE estimator, which is \cite{Kay-Est.Theory}
\begin{equation}
\hat{\ve{x}}_{\text{L}}  = \m{C}_{\ve{x}\ve{x}}\m{H}^H(\m{H}\m{C}_{\ve{x}\ve{x}}\m{H}^H+\m{C}_{\ve{n}\ve{n}})^{-1}\ve{y} =  \m{E}_{\text{L}}\ve{y}.  \label{equ:CWCULMMSE002c} 
\end{equation}
Let $\ve{e}_{\text{L},i}^H\in \mathbb{C}^{1 \times m}$ be the $i^{th}$ row of $\m{E}_{\text{L}}$, then the conditional mean and variance are given by \eqref{equ:CWCULMMSE016a} and \eqref{equ:CWCULMMSE102}, respectively, where $\ve{e}_{\text{L},i}^H$ has to be inserted for $\ve{e}_i^H$. A known property of the LMMSE estimator is that $\alpha_{\text{L},i} = \ve{e}_{\text{L},i}^H\ve{h}_i$ is real valued, and  in general smaller than $1$. Hence, $\hat{x}_{\text{L},i}$ is conditionally biased according to \eqref{equ:CWCULMMSE016a}.

We now turn to the \underline{C}WCU \underline{L}MMSE estimator which is given by \cite{Huemer_2014_Asilomar}, \cite{Lang_2015_Eurocast}
\begin{equation}
\hat{\ve{x}}_{\text{CL}}  = \m{D}\m{C}_{\ve{x}\ve{x}}\m{H}^H(\m{H}\m{C}_{\ve{x}\ve{x}}\m{H}^H+\m{C}_{\ve{n}\ve{n}})^{-1}\ve{y} =  \m{E}_{\text{CL}}\ve{y}, \label{equ:CWCULMMSE002b}
\end{equation}
where the elements of the real diagonal matrix $\m{D}$ are $[\m{D}]_{i,i} = 1/\alpha_{\text{L},i}$.
The CWCU LMMSE estimator in \eqref{equ:CWCULMMSE002b} and the LMMSE estimator in \eqref{equ:CWCULMMSE002c}  are connected via
\begin{equation}
\hat{\ve{x}}_{\text{CL}} = \m{D} \m{E}_{\text{L}}\ve{y} = \m{D} \hat{\ve{x}}_{\text{L}}.   \label{equ:CWCULMMSE002ba}
\end{equation}
Let $\ve{e}_{\text{CL},i}^H\in \mathbb{C}^{1 \times m}$ be the $i^{th}$ row of $\m{E}_{\text{CL}}$, then it holds that $\ve{e}_{\text{L},i}^H = \alpha_{\text{L},i} \ve{e}_{\text{CL},i}^H$, $\hat{x}_{\text{L},i} = \alpha_{\text{L},i} \hat{x}_{\text{CL},i}$ and $\text{var}(\hat{x}_{\text{L},i}|x_i) = \alpha_{\text{L},i}^2 \text{var}(\hat{x}_{\text{CL},i}|x_i)$.

In contrast to the LMMSE estimator, the CWCU LMMSE estimator fulfills $\ve{e}_{\text{CL},i}^H\ve{h}_i=1$. This property makes \eqref{equ:CWCULMMSE016a} equal to $E_{\ve{y}|x_i}[\hat{x}_i|x_i] = x_i$ (which is the CWCU constraint in \eqref{equ:cuLMMSE047}).  Hence, $\hat{x}_{\text{CL},i}$ is conditionally unbiased. The conditional mean and variance of the CWCU LMMSE estimator are given by \eqref{equ:CWCULMMSE016a} and \eqref{equ:CWCULMMSE102}, respectively, where $\ve{e}_{\text{CL},i}^H$ has to be inserted for $\ve{e}_i^H$. Inserting these conditional properties into \eqref{equ:CWCULMMSE105} yields
\begin{align}
 &p(\hat{x}_{\text{CL},i}|s^{(q)})  \nonumber \\
 \hspace{2mm} &= \frac{1}{ \pi \text{var}(\hat{x}_{\text{CL},i}|s^{(q)})} \text{e}^{-\frac{1}{ \text{var}(\hat{x}_{\text{CL},i}|s^{(q)})} \left| \hat{x}_{\text{CL},i} -  E[\hat{x}_{\text{CL},i}|s^{(q)}] \right|^2} \nonumber \\
 & \hspace{2mm}= \frac{\alpha_{\text{L},i}^2}{ \pi \text{var}(\hat{x}_{\text{L},i}|s^{(q)})} \text{e}^{-\frac{\alpha_{\text{L},i}^2}{ \text{var}(\hat{x}_{\text{L},i}|s^{(q)})} \big| \alpha_{\text{L},i}^{-1}\big(\hat{x}_{\text{L},i} -  \alpha_{\text{L},i} s^{(q)} \big) \big|^2} \nonumber \\ 
 & \hspace{2mm}= \frac{\alpha_{\text{L},i}^2}{ \pi \text{var}(\hat{x}_{\text{L},i}|s^{(q)})} \text{e}^{-\frac{1}{ \text{var}(\hat{x}_{\text{L},i}|s^{(q)})} \left| \hat{x}_{\text{L},i} - E[\hat{x}_{\text{L},i}|s^{(q)}] \right|^2} \nonumber \\ 
 & \hspace{2mm}= \alpha_{\text{L},i}^2 p(\hat{x}_{\text{L},i}|s^{(q)}), \label{equ:CWCULMMSE105a}
\end{align}
which holds for any symbol $s^{(q)}$. Hence, for a given $\ve{y}$ the probability density $p(\hat{x}_{\text{CL},i}|s^{(q)})$ of the CWCU LMMSE estimator and $p(\hat{x}_{\text{L},i}|s^{(q)})$ of the LMMSE estimator for any $s^{(q)}$ only differ by the constant scaling factor $\alpha_{\text{L},i}^2$. This constant scaling factor does not depend on the symbol $s^{(q)}$ and it appears in the numerator and the denominator of \eqref{equ:LLR001}, thus cancelling out. Hence, the LLRs of the CWCU LMMSE estimates and the LMMSE estimates are equal for proper constellation diagrams. As a consequence, the resulting bit error ratios (BERs) of the LMMSE and the CWCU LMMSE estimators are also the same, although the BMSE of the LMMSE estimator is in general lower than that of the CWCU LMMSE estimator.

\section{Widely Linear Estimation of Improper Data} \label{sec:Widely Linear Estimation of Improper Data}

We now turn to improper constellation diagrams such as 8-QAM. 
In such scenarios it is advantageous to use widely linear estimators, which can incorporate information about the improperness of the data.
A general widely linear estimator in augmented notation is 
\begin{equation}
	\underline{\hat{\ve{x}}} = \begin{bmatrix} \hat{\ve{x}} \\ \hat{\ve{x}}^{*} \end{bmatrix} = \begin{bmatrix}  \m{E} & \m{F} \\ \m{F}^{*} & \m{E}^{*} \end{bmatrix}  \begin{bmatrix} \ve{y} \\ \ve{y}^{*} \end{bmatrix}  = \underline{\m{E}} \underline{\ve{y}}, \label{equ:CWCU_Journal106}
\end{equation}
where $(\cdot)^{*}$ denotes the complex conjugate. For an introduction to the augmented form and widely linear estimation we refer to \cite{Schreier-2011}. Isolating the $i^{th}$ element of \eqref{equ:CWCU_Journal106} yields $\hat{x}_{i} =  \ve{e}_i^H  \underline{\ve{y}}$,
where $\ve{e}_i^H\in \mathbb{C}^{1 \times 2m}$ is the $i^{th}$ row of $\underline{\m{E}}$. The augmented version is given by
\begin{equation}
\underline{\hat{\ve{x}}}_{i} = \begin{bmatrix} \hat{x}_i \\ \hat{x}_i^{*} \end{bmatrix} =\begin{bmatrix}
 \ve{e}_i^H \\ \ve{e}_{i+n}^H
\end{bmatrix}\underline{\ve{y}} = \underline{\m{E}}_i^H \underline{\ve{y}},  \label{equ:LLR013}
\end{equation}
where the rows of $\underline{\m{E}}_i^H$ are given by the $i^{th}$ and the $(i+n)^{th}$ row of the augmented estimator matrix $\underline{\m{E}}$. The augmented version of \eqref{equ:CWCULMMSE015a} is
\begin{equation}
\underline{\ve{y}} =\begin{bmatrix} \m{H} & \m{0} \\  \m{0} & \m{H}^* \end{bmatrix} \underline{\ve{x}} + \underline{\ve{n}} = \underline{\m{H}}\underline{\ve{x}} + \underline{\ve{n}} = \underline{\ve{H}}_i \underline{\ve{x}}_i + \underline{\bar{\m{H}}}_i\underline{\bar{\ve{x}}}_i + \underline{\ve{n}}, \label{equ:LLR018}
\end{equation}
where
\begin{align}
\underline{\ve{H}}_i =& \begin{bmatrix}
\ve{h}_i & \hspace{-2mm} \ve{0} \\ \ve{0} & \hspace{-2mm}\ve{h}_i^*
\end{bmatrix}, & \hspace{-2mm}
\underline{\ve{x}}_i =& \begin{bmatrix}
x_i \\ x_i^*
\end{bmatrix}, & \hspace{-2mm} 
\underline{\bar{\m{H}}}_i =& \begin{bmatrix}
\bar{\m{H}}_i & \hspace{-2mm}\m{0} \\
\m{0} &\hspace{-2mm} \bar{\m{H}}_i^* 
\end{bmatrix}, & \hspace{-2mm}
\underline{\bar{\ve{x}}}_i =& \begin{bmatrix}
\bar{\ve{x}}_i  \\  \bar{\ve{x}}_i^*
\end{bmatrix}. \nonumber
\end{align}
With \eqref{equ:LLR018}, \eqref{equ:LLR013} can be rewritten according to
\begin{equation}
\underline{\hat{\ve{x}}}_{i} = \underline{\m{E}}_i^H \underline{\ve{H}}_i \underline{\ve{x}}_i + \underline{\m{E}}_i^H\underline{\bar{\m{H}}}_i\underline{\bar{\ve{x}}}_i + \underline{\m{E}}_i^H \underline{\ve{n}}.   \label{equ:LLR020}
\end{equation}
For zero mean and statistically independent elements of $\ve{x}$, the conditional augmented expected vector of $\underline{\hat{\ve{x}}}_i$ follows to
\begin{equation}
E[\underline{\hat{\ve{x}}}_{i} | x_i] =  \underline{\m{E}}_i^H \underline{\ve{H}}_i \underline{\ve{x}}_i = \underline{\m{\alpha}}_{i} \underline{\ve{x}}_i. \label{equ:LLR021}
\end{equation}
From \eqref{equ:LLR020} and \eqref{equ:LLR021}, the augmented conditional covariance matrix of $\hat{x}_{i}$ is
\begin{align}
\underline{\m{C}}_{\hat{x}_{i}\hat{x}_{i}|x_i} &=E[ (\underline{\hat{\ve{x}}}_{i} - E[\underline{\hat{\ve{x}}}_{i} | x_i] )  (\underline{\hat{\ve{x}}}_{i} - E[\underline{\hat{\ve{x}}}_{i} | x_i] )^H  |x_i ] \nonumber \\
& = E[ \underline{\m{E}}_i^H \left( \underline{\bar{\m{H}}}_i \underline{\bar{\ve{x}}}_i  +  \underline{\ve{n}} \right)   \left( \underline{\bar{\m{H}}}_i \underline{\bar{\ve{x}}}_i +  \underline{\ve{n}} \right)^H \underline{\m{E}}_i  |x_i ] \nonumber \\
&=  \underline{\m{E}}_i^H  \left( \underline{\bar{\m{H}}}_i \underline{\m{C}}_{\bar{\ve{x}}_i\bar{\ve{x}}_i} \underline{\bar{\m{H}}}_i^H + \underline{\m{C}}_{\ve{n}\ve{n}}   \right) \underline{\m{E}}_i. \label{equ:LLR022}
\end{align}
Similar to the linear case in \eqref{equ:CWCULMMSE102}, \eqref{equ:LLR022} is independent of $x_i$. Particular realizations for \eqref{equ:LLR021} and \eqref{equ:LLR022} can be obtained by inserting $ \underline{\m{E}}_i^H$ of a concrete estimator. Such a particular estimator could be the \underline{WL}MMSE estimator, whose augmented form is \cite{Schreier-2011}
\begin{equation}
\underline{\hat{\ve{x}}}_{\text{WL}} = \underline{\m{C}}_{\ve{x}\ve{y}}\underline{\m{C}}_{\ve{y}\ve{y}}^{-1}  \underline{\ve{y}} = \underline{\m{E}}_{\text{WL}}  \underline{\ve{y}}. \label{equ:LLR026}
\end{equation}
Then, the augmented $i^{th}$ estimate is given by
\begin{equation}
  \underline{\hat{\ve{x}}}_{\mathrm{WL},i} = \underline{\m{C}}_{x_i\ve{y}}\underline{\m{C}}_{\ve{y}\ve{y}}^{-1}  \underline{\ve{y}} =
  \underline{\m{E}}_{\mathrm{WL},i}^H \underline{\ve{y}},\label{equ:CWCU_Journal028d}  
\end{equation}
where the rows of $\underline{\m{E}}_{\mathrm{WL},i}^H$ are the $i^{th}$ and the $(i+n)^{th}$ row of $\underline{\m{E}}_{\text{WL}}$ in \eqref{equ:LLR026}. 
For the WLMMSE estimator, $\underline{\m{\alpha}}_{\mathrm{WL},i} = \underline{\m{E}}_{\mathrm{WL},i}^H \underline{\ve{H}}_i$ is in general not equal to the identity matrix. Hence, according to \eqref{equ:LLR021}, $\underline{\hat{\ve{x}}}_{\mathrm{WL},i}$ is conditionally biased. 

We now turn to the \underline{C}WCU \underline{WL}MMSE estimator, whose augmented $i^{th}$ estimate is \cite{CWCU_WLMMSE_Journal_Huemer_2016}
\begin{equation}
  \underline{\hat{\ve{x}}}_{\mathrm{CWL},i} = \underline{\m{C}}_{x_i x_i}\left(  \underline{\m{C}}_{x_i\ve{y}}\underline{\m{C}}_{\ve{y}\ve{y}}^{-1}\underline{\m{C}}_{\ve{y}x_i}\right)^{-1}\underline{\m{C}}_{x_i\ve{y}}\underline{\m{C}}_{\ve{y}\ve{y}}^{-1}  \underline{\ve{y}}.\label{equ:CWCU_Journal028c}  
\end{equation}
For statistically independent elements of $\ve{x}$, it holds that $\underline{\m{C}}_{\ve{y}x_i}=\underline{\ve{H}}_i \underline{\m{C}}_{x_i x_i}$ and \eqref{equ:CWCU_Journal028c} can be reformulated as
\begin{align}
& \underline{\hat{\ve{x}}}_{\mathrm{CWL},i} = \underline{\m{C}}_{x_i x_i}\bigg(  \underbrace{\underline{\m{C}}_{x_i\ve{y}}\underline{\m{C}}_{\ve{y}\ve{y}}^{-1}}_{\underline{\m{E}}_{\mathrm{WL},i}^H}\underline{\m{C}}_{\ve{y}x_i}\bigg)^{-1}\underbrace{\underline{\m{C}}_{x_i\ve{y}}\underline{\m{C}}_{\ve{y}\ve{y}}^{-1}}_{\underline{\m{E}}_{\mathrm{WL},i}^H} \underline{\ve{y}} \nonumber \\
&\hspace{1mm}= \left(  \underline{\m{E}}_{\mathrm{WL},i}^H\underline{\m{C}}_{\ve{y}x_i} \underline{\m{C}}_{x_i x_i}^{-1}\right)^{-1}\underline{\m{E}}_{\mathrm{WL},i}^H \underline{\ve{y}} \nonumber\\
&\hspace{1mm}=  \big( \underbrace{ \underline{\m{E}}_{\mathrm{WL},i}^H\underline{\ve{H}}_i}_{\underline{\m{\alpha}}_{\mathrm{WL},i}} \big)^{-1}\underline{\m{E}}_{\mathrm{WL},i}^H \underline{\ve{y}} 
= \underline{\m{\alpha}}_{\mathrm{WL},i}^{-1}\underline{\hat{\ve{x}}}_{\mathrm{WL},i}. \label{equ:CWCU_Journal028e}    
\end{align}
Similar to the linear case in \eqref{equ:CWCULMMSE002ba}, the CWCU WLMMSE estimator is determined by the WLMMSE estimator times a term that corrects for the conditional bias. It follows from \eqref{equ:CWCU_Journal028e}, that the augmented conditional covariance matrix of the CWCU WLMMSE estimator can be derived from the one of the WLMMSE estimator according to
\begin{align}
\underline{\m{C}}_{\hat{x}_{i}\hat{x}_{i}|x_i, \text{CWL}} = \underline{\m{\alpha}}_{\mathrm{WL},i}^{-1} \underline{\m{C}}_{\hat{x}_{i}\hat{x}_{i}|x_i, \text{WL}} \left(\underline{\m{\alpha}}_{\mathrm{WL},i}^H\right)^{-1}. \label{equ:CWCU_Journal028f}
\end{align}
With these conditional properties, it is possible to evaluate the LLRs by utilizing the general complex Gaussian density function
\begin{align}
&p(\hat{x}_i|s^{(q)}) =   \frac{1}{\sqrt{\pi^2 \text{det}(\underline{\m{C}}_{\hat{x}_{i}\hat{x}_{i}|s^{(q)}})}} \nonumber \\
& \hspace{10mm} \cdot \text{e}^{-\frac{1}{2} \left( \underline{\hat{\ve{x}}}_{i} - E[\underline{\hat{\ve{x}}}_{i} | s^{(q)}] \right) ^H \underline{\m{C}}_{\hat{x}_{i}\hat{x}_{i}|s^{(q)}}^{-1}  \left( \underline{\hat{\ve{x}}}_{i} - E[\underline{\hat{\ve{x}}}_{i} | s^{(q)}] \right) }.  \label{equ:LLR025}
\end{align}
In analogy to the linear case in \eqref{equ:CWCULMMSE105a} it will now be shown that $p(\hat{x}_{\text{WL},i}|s^{(q)})$ of the WLMMSE estimator and $p(\hat{x}_{\text{CWL},i}|s^{(q)})$ of the CWCU WLMMSE estimator only differ by a constant factor. By utilizing \eqref{equ:CWCU_Journal028e}  and \eqref{equ:CWCU_Journal028f} the exponent of \eqref{equ:LLR025} for the CWCU WLMMSE estimator can be rearranged to
\begin{align}
-\frac{1}{2}& \left( \underline{\hat{\ve{x}}}_{\mathrm{CWL},i} - E[\underline{\hat{\ve{x}}}_{\mathrm{CWL},i} | s^{(q)}] \right) ^H \underline{\m{C}}_{\hat{x}_{i}\hat{x}_{i}|s^{(q)}, \text{CWL}}^{-1}  \nonumber \\
& \hspace{5mm} \cdot \left( \underline{\hat{\ve{x}}}_{\mathrm{CWL},i} - E[\underline{\hat{\ve{x}}}_{\mathrm{CWL},i} | s^{(q)}] \right) \nonumber \\
&=-\frac{1}{2} \left( \underline{\hat{\ve{x}}}_{\mathrm{WL},i} - E[\underline{\hat{\ve{x}}}_{\mathrm{WL},i} | s^{(q)}] \right) ^H \left(\underline{\m{\alpha}}_{\mathrm{WL},i}^H\right)^{-1} \underline{\m{\alpha}}_{\mathrm{WL},i}^{H} \nonumber  \\
& \hspace{5mm} \cdot \underline{\m{C}}_{\hat{x}_{i}\hat{x}_{i}|s^{(q)}, \text{WL}}^{-1}\underline{\m{\alpha}}_{\mathrm{WL},i} \underline{\m{\alpha}}_{\mathrm{WL},i}^{-1} \left( \underline{\hat{\ve{x}}}_{\mathrm{WL},i} - E[\underline{\hat{\ve{x}}}_{\mathrm{WL},i} | s^{(q)}] \right)\nonumber \\
&=-\frac{1}{2} \left(\ \underline{\hat{\ve{x}}}_{\mathrm{WL},i} -  E[\underline{\hat{\ve{x}}}_{\mathrm{WL},i} |s^{(q)}] \right) ^H \underline{\m{C}}_{\hat{x}_{i}\hat{x}_{i}|s^{(q)}, \text{WL}}^{-1}  \nonumber \\
& \hspace{5mm} \cdot \left(  \underline{\hat{\ve{x}}}_{\mathrm{WL},i} - E[\underline{\hat{\ve{x}}}_{\mathrm{WL},i} |s^{(q)}]  \right). \label{equ:LLR028}
\end{align}
This result shows that the exponent of \eqref{equ:LLR025} is identical for the CWCU WLMMSE estimator and the WLMMSE estimator for a given $\ve{y}$. The prefactor of \eqref{equ:LLR025} for the CWCU WLMMSE estimator follows to
\begin{align}
 &\frac{1}{\sqrt{\pi^2 \text{det}(\underline{\m{C}}_{\hat{x}_{i}\hat{x}_{i}|s^{(q)}, \text{CWL}})}} \nonumber \\
 &\hspace{5mm}= \frac{1}{\sqrt{\pi^2 \text{det}\left(\underline{\m{\alpha}}_{\mathrm{WL},i}^{-1} \underline{\m{C}}_{\hat{x}_{i}\hat{x}_{i}|s^{(q)}, \text{WL}} \left(\underline{\m{\alpha}}_{\mathrm{WL},i}^H\right)^{-1}\right)}} \nonumber \\
 &\hspace{5mm}=\frac{1}{\sqrt{\pi^2  \text{det}( \underline{\m{C}}_{\hat{x}_{i}\hat{x}_{i}|s^{(q)}, \text{WL}}) |\text{det}(\underline{\m{\alpha}}_{\mathrm{WL},i}^{-1})|^2}} \nonumber\\
 &\hspace{5mm}=\frac{\left|\text{det}(\underline{\m{\alpha}}_{\mathrm{WL},i})\right|}{ \sqrt{\pi^2  \text{det}( \underline{\m{C}}_{\hat{x}_{i}\hat{x}_{i}|s^{(q)}, \text{WL}}) }}. \label{equ:LLR027}
\end{align}
Like in the linear case in \eqref{equ:CWCULMMSE105a}, the prefactors of the CWCU WLMMSE estimator and the WLMMSE estimator only differ by a constant real factor. 
This factor does not depend on the symbol $s^{(q)}$ and it appears in the numerator and the denominator of \eqref{equ:LLR001}, thus cancelling out in the determination of the LLRs. This leads to the result that the LLRs derived from the CWCU WLMMSE estimates and the WLMMSE estimates are exactly the same. Although the WLMMSE estimator in general features a lower BMSE, the BER performance of the WLMMSE and the CWCU WLMMSE estimator are identical.

\section{Simulation Example} \label{sec:Propriety of the estimated data symbols}

We give a simulation example were we use the unique word orthogonal frequency division multiplexing (UW-OFDM) framework described in \cite{UW_OFDM}, \cite{UW_OFDM_2}. Like classical OFDM, UW-OFDM is a block based transmission scheme where in our particular setup at the receive side a data vector $\ve{d} \in \mathbb{C}^{36 \times 1}$ is estimated based on a received block $\vef{y} \in \mathbb{C}^{52 \times 1}$ of frequency domain samples. We choose UW-OFDM since the estimator matrices are in general full matrices instead of diagonal matrices as in classical OFDM, such that the problem can be considered a more demanding and general one compared to the data estimation problem in classical OFDM systems. Hence, this framework is well suited for studying general effects of CWCU estimators. The system model for the transmission of one data block is given by 
\begin{equation}
\vef{y} = \mf{H}\m{G} \ve{d} + \vef{v},
\end{equation}
where $\mf{H} \in \mathbb{C}^{52 \times 52}$ is the diagonal channel matrix. $\m{G} \in \mathbb{C}^{52 \times 36}$ is a so called generator matrix, for details cf. \cite{UW_OFDM}, \cite{UW_OFDM_2}, $\ve{d}$ is a vector of improper 8-QAM symbols and $\vef{v}$ is a frequency domain noise vector. Note that every assumption made in Sec.~\ref{sec:Linear estimation of proper data} and Sec.~\ref{sec:Widely Linear Estimation of Improper Data} holds in this example: The data and the measurements are connected via a linear model, and the Gaussian assumption of $p(\hat{x}_i|s^{(q)})$ is valid due to central limit theorem arguments (note that the data vector length is 36 in this example).

\begin{figure}[!t]
\centering
\includegraphics[width=3.35in]{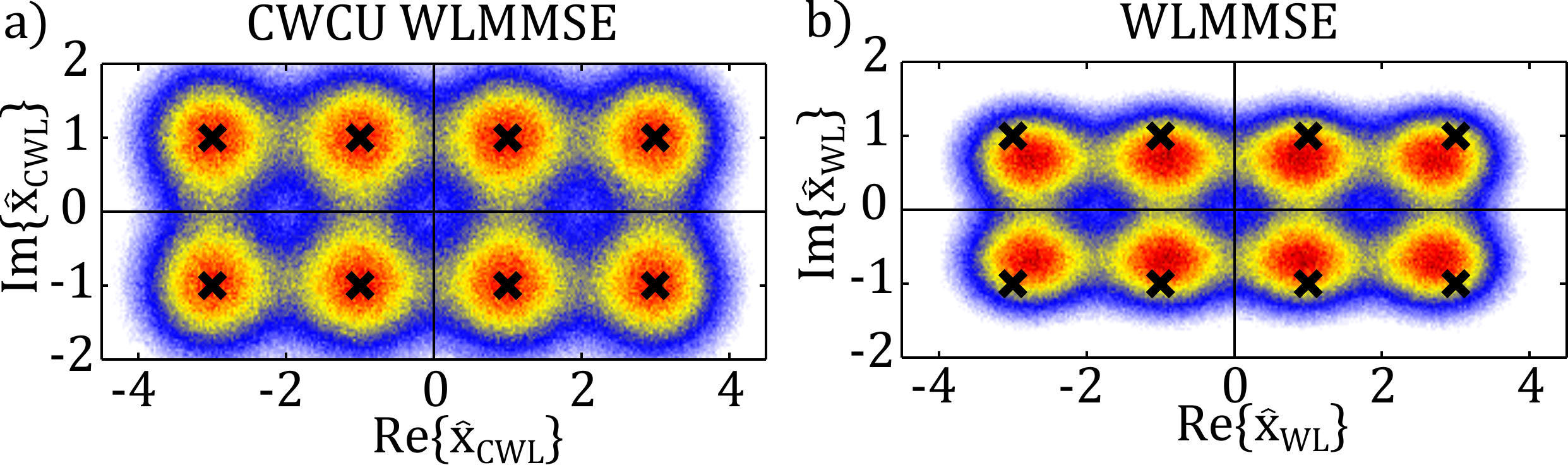}
\caption{Relative frequencies of the CWCU WLMMSE estimates in (a), and the WLMMSE estimates in (b). The black crosses mark the original 8-QAM constellation points.}
\label{fig:PDFs}
\end{figure}

In the simulation, UW-OFDM symbols are transmitted over an AWGN channel $\mf{H} = \m{I}$ and further processed by the WLMMSE estimator and the CWCU WLMMSE estimator, respectively. These estimators feature different properties of the estimated data symbols. According to \cite{Schreier-2011}, the estimates conditioned on a given $s^{(q)}$ are proper, if the off-diagonal elements of $\underline{\m{C}}_{\hat{x}_{i}\hat{x}_{i}|s^{(q)}}$ are zero, which holds true for 8-QAM symbols transmitted over the AWGN channel and received by the CWCU WLMMSE estimator. The corresponding relative frequencies of $\underline{\hat{\ve{x}}}_{\mathrm{CWL},i}$ are shown in Fig.~\ref{fig:PDFs}a. One can see that the estimates are centered around the true constellation points since the CWCU WLMMSE estimator fulfills the CWCU constraints in \eqref{equ:cuLMMSE047}. Furthermore, the estimates conditioned on a specific transmit symbol are properly distributed. In Fig.~\ref{fig:PDFs}b, the relative frequencies of the WLMMSE estimates are shown. In contrast to the CWCU WLMMSE estimates, the WLMMSE estimates conditioned on a specific transmit symbol are neither centered around the true constellation points nor are they properly distributed. However, due to the close connection between the CWCU WLMMSE estimator and the WLMMSE estimator, the resulting LLRs are identical as shown above. Moreover, since the CWCU WLMMSE estimates for a given $s^{(q)}$ are proper, it is sufficient to use the proper complex Gaussian PDF in \eqref{equ:CWCULMMSE105} instead of the general Gaussian PDF in \eqref{equ:LLR025} as basis for the LLR determination. 
Hence, the LLR determination of the CWCU WLMMSE estimates is less computationally demanding than for the WLMMSE estimates without any loss in BER performance. 

In \cite{CWCU_WLMMSE_Journal_Huemer_2016}, we also confirmed via simulation, that for frequency selective channels the CWCU WLMMSE estimates conditioned on a given transmit symbol $s^{(q)}$ are practically proper again for all investigated channel realizations. The off-diagonal elements of $\underline{\m{C}}_{\hat{x}_{i}\hat{x}_{i}|s^{(q)}}$ are smaller than the main diagonal elements by at least a factor of $10^{3}$. Rounding the off-diagonal elements to zero and applying the proper complex Gaussian PDF in \eqref{equ:CWCULMMSE105} instead of the general Gaussian PDF in \eqref{equ:LLR025} as basis for the LLR determination leads to a BER without any noticeable loss in performance.

\section{Conclusion}

In this paper, we proved that the CWCU LMMSE estimates result in the same LLRs as the LMMSE estimates for proper constellation diagrams such as QPSK or 16-QAM. As a consequence, the resulting BER performance of the CWCU LMMSE estimator and the LMMSE estimator is also the same, even though the two estimators fulfill different unbiased constrains and yield a different BMSE. For improper constellation diagrams such as 8-QAM, we showed that the same statements also hold for the relationships between the widely linear counterparts, the CWCU WLMMSE and WLMMSE estimators. A simulation example was presented, revealing different statistical properties of WLMMSE and CWCU WLMMSE data estimates. An interesting outcome is that the CWCU WLMMSE estimator offers a complexity advantage in the LLR determination over the WLMMSE estimator without a loss in BER performance.



\end{document}